\documentclass{tran-l}
\usepackage[dvips]{graphicx}
\usepackage{amscd}
\usepackage{amsmath}
\usepackage{amsxtra}
\usepackage{amsfonts}
\usepackage{amssymb}
\newtheorem{theorem}{Theorem}[section]
\newtheorem{corollary}[theorem]{Corollary}
\newtheorem{lemma}[theorem]{Lemma}
\newtheorem{proposition}[theorem]{Proposition}
\theoremstyle{definition}
\newtheorem{definition}[theorem]{Definition}
\newtheorem{remark}[theorem]{Remark}

\newtheorem{example}[theorem]{Example}
\theoremstyle{remark}

\renewcommand{\theclaim}{\textup{\theclaim}}

\numberwithin{equation}{section}

\def\openone%{\hbox{\upshape \small1\kern-3.3pt\normalsize1}}

{\mathchoice

{\hbox{\upshape \small1\kern-3.3pt\normalsize1}}

{\hbox{\upshape \small1\kern-3.3pt\normalsize1}}

{\hbox{\upshape \tiny1\kern-2.3pt\SMALL1}}

{\hbox{\upshape \Tiny1\kern-2pt\tiny1}}}

\makeatletter

\newbox\ipbox

\newcommand{\ip}[2]{\left\langle #1\,|\,#2\right\rangle}

\newcommand{\diracb}[1]{\left\langle #1\mathrel{\mathchoice

{\setbox\ipbox=\hbox{$\displaystyle \left\langle\mathstrut #1\right.$}

\vrule height\ht\ipbox width0.25pt depth\dp\ipbox}

{\setbox\ipbox=\hbox{$\textstyle \left\langle\mathstrut #1\right.$}

\vrule height\ht\ipbox width0.25pt depth\dp\ipbox}

{\setbox\ipbox=\hbox{$\scriptstyle \left\langle\mathstrut #1\right.$}

\vrule height\ht\ipbox width0.25pt depth\dp\ipbox}

{\setbox\ipbox=\hbox{$\scriptscriptstyle \left\langle\mathstrut #1\right.$}

\vrule height\ht\ipbox width0.25pt depth\dp\ipbox}

}\right. }

\newcommand{\dirack}[1]{\left. \mathrel{\mathchoice

{\setbox\ipbox=\hbox{$\displaystyle \left.\mathstrut #1\right\rangle$}

\vrule height\ht\ipbox width0.25pt depth\dp\ipbox}

{\setbox\ipbox=\hbox{$\textstyle \left.\mathstrut #1\right\rangle$}

\vrule height\ht\ipbox width0.25pt depth\dp\ipbox}

{\setbox\ipbox=\hbox{$\scriptstyle \left.\mathstrut #1\right\rangle$}

\vrule height\ht\ipbox width0.25pt depth\dp\ipbox}

{\setbox\ipbox=\hbox{$\scriptscriptstyle \left.\mathstrut #1\right\rangle$}

\vrule height\ht\ipbox width0.25pt depth\dp\ipbox}

} #1\right\rangle}

\newcommand{\lonet}{L^{1}\left(  \mathbb{T}\right)}

\newcommand{\linft}{L^{\infty}\left(  \mathbb{T}\right)}

\newcommand{\ltwor}{L^{2}\left(\mathbb{R}\right)}

\newcommand{\ltwoz}{l^2\left(\mathbb{Z}\right)}

\newcommand{\cj}[1]{\overline{#1}}

\newcommand{\Zu}{\mathbb{Z}}

\newcommand{\SN}{\mathcal{S}_N}
\newcommand{\ltwoez}{L^2(E,l^2(\mathbb{Z}))}
\newcommand{\bz}{\mathbb{Z}}
\newcommand{\br}{\mathbb{R}}
\newcommand{\bc}{\mathbb{C}}
\newcommand{\bt}{\mathbb{T}}

\input cyracc.def

\begin{document}
\title[The Baumslag Solitar group]{Low-pass filters and representations of the Baumslag Solitar group}
\author{Dorin Ervin Dutkay}
\address{ Department of Mathematics\\
Hill Center-Busch Campus\\
Rutgers, The State University of New Jersey\\
110 Frelinghuysen Rd\\
Piscataway, NJ 08854-8019, USA } \email{Dorin Ervin Dutkay: ddutkay@math.rutgers.edu}
\thanks{This work was supported in part bt NSF grant DMS0457491}

\renewcommand{\subjclassname}{\textup{2000} Mathematics Subject Classification}

\begin{abstract}
We analyze representations of the Baumslag Solitar group
$$BS(1,N)=\langle u,t\,|\,utu^{-1}=t^N\rangle$$ that admit wavelets and show how such representations can be constructed from a given low-pass filter.
We describe the direct integral decomposition for some examples
and derive from it a general criterion for the existence of
solutions for scaling equations. As another application, we construct a Fourier
transform for some Hausdorff measures.
\end{abstract}

\maketitle \tableofcontents
\section{\label{Intro}Introduction and notations}
Wavelet theory was developed initially for $\ltwor$ to construct
orthonormal bases that have good localization properties (see
\cite{Dau92} for more details). In $\ltwor$ we have two operators:
the translation operator
$$Tf(x)=f(x-1),\quad(x\in\br,f\in\ltwor),$$
and the dilation operator
$$Uf(x)=\frac{1}{\sqrt{N}}f\left(\frac{x}{N}\right),\quad(x\in\br,f\in\ltwor),$$
where $N\geq 2$ is an integer called the scale.
\par
A wavelet is a finite set $\{\psi_1,...,\psi_p\}$ of functions in
$\ltwor$ such that
$$\{U^jT^k\psi_i\,|\,j,k\in\bz,i\in\{1,...,p\}\}$$
is an orthonormal basis for $\ltwor$.
\par
The operators satisfy the commutation relation $UTU^{-1}=T^N$,
hence we are dealing with a representation of the Baumslag Solitar
group $BS(1,N)$. Here is the definition and some elementary facts:
\par

 The Baumslag Solitar group $BS(1,N)$ is the group with two
generators $u$ and $t$ and one relation
$$utu^{-1}=t^N.$$
It can also be described the semidirect product of the group of
$N$-adic numbers $\mathbb{Z}[1/N]$, with $\mathbb{Z}$, the action
of $\mathbb{Z}$ on $\mathbb{Z}[1/N]$ being given by multiplication
by $N$:
$$\alpha_i\left(\frac{k}{N^p}\right)=N^i\frac{k}{N^p},\quad(i\in\mathbb{Z},\frac{k}{N^p}\in\mathbb{Z}[1/N]).$$
The elements in $\mathbb{Z}[1/N]$ of the form $k/N^p$ correspond
to $u^pt^ku^{-p}$ which we denote by $t_{k/N^p}$, and the elements
in $j\in\mathbb{Z}$, correspond to $u^j$. Each element in the
group can be written uniquely as $u^jt_d$ with $j\in\mathbb{Z}$
and $d\in\mathbb{Z}[1/N]$. The multiplication rule is
$$(u^jt_d)(u^{j'}t_{d'})=u^{j+j'}t_{N^{-j'}d+d'}.$$
The Baumslag Solitar group $BS(1,N)$ can be regarded also as an
$ax+b$-group with $a$ an integer power of $N$ and $b$ $N$-adic.
\par
The action $\alpha$ is given by conjugation by $u$:
$$\alpha_i(t_d)=u^it_du^{-i},\quad(d\in\mathbb{Z},i\in\mathbb{Z}).$$

\par
The dual of the group $\mathbb{Z}[1/N]$ is the solenoid:
$$\mathcal{S}_N:=\{(z_n)_{n\geq
0}\,|\,z_{n+1}^N=z_n,|z_n|=1\mbox{ for }n\geq 0\}.$$
Sometimes it will be useful to index the components by
$\mathbb{Z}$ and identify $\mathcal{S}_N$ with
$$\mathcal{S}_N:=\{(z_n)_{n\in\mathbb{Z}}\,|\,z_{n+1}^N=z_n,|z_n|=1\mbox{ for }n\in\bz\},$$ the
identification being given by
$$(z_n)_{n\geq 0}\leftrightarrow (z_n)_{n\in\mathbb{Z}},\mbox{ with }z_{-n}=z_0^{N^n}\mbox {for }n>0.$$
The duality is given by
\begin{equation}\label{eqdual}
\ip{\frac{k}{N^p}}{(z_n)_{n\in\mathbb{Z}}}=z_p^k,\quad(\frac{k}{N^p}\in\mathbb{Z}[1/N],(z_n)_{n\in\Zu}\in\SN).
\end{equation}
The dual action, $\hat\alpha$ of $\mathbb{Z}$ on $\mathcal{S}_N$,
is given by the shift $S$ on $\SN$:
$$S(z_n)_{n\in\Zu}=(z_{n-1})_{n\in\Zu},\quad
\hat\alpha_i(z_n)_{n\in\Zu}=S^i(z_n)_{n\in\Zu},\quad((z_n)_{n\in\Zu}\in\SN).$$
The dual of the inclusion $i:\bz[1/N]\rightarrow\br$ is $\hat
i:\br\rightarrow\SN$,
$$\hat i(x)=(e^{-2\pi\mbox{i}N^{-n}x})_{n\in\bz},\quad(x\in\br).$$
\par
If $\pi$ is a unitary representation of the group $BS(1,N)$ on
some Hilbert space we will use capital letters to denote the
corresponding operators: $\pi(u)=:U$, $\pi(t)=:T$, $\pi(t_d)=T_d$
for $d\in\bz[1/N]$.
\par
Representations of the Baumslag Solitar group that admit wavelet
bases can be constructed also on some other spaces such as:
$\ltwor\oplus...\oplus\ltwor$ (\cite{BDP}, \cite{Dut2})) or some
fractal spaces (\cite{DutJo}).
\par
We define here the concepts of wavelet theory in the abstract
setting, when we have a representation of the group $BS(1,N)$ on
some Hilbert space.
\par

Let $\pi$ be a representation of the group $BS(1,N)$ on some
Hilbert space $H$. A wavelet for the representation $\pi$ is a
finite set $\{\psi_1,...,\psi_n\}$ of vectors in $H$ such that
$$\{\pi(u^jt^k)\psi_i\,|\,j,k\in\bz,i\in\{1,...,n\}\}$$
is an orthonormal basis for $H$.
\par
We will see in section \ref{rep} that the representations that
admit wavelets are faithful and weakly equivalent to the right
regular representation of the group.
\par
The main technique to construct wavelets is by a {\em
multiresolution analysis}, which is a sequence of subspaces
$(V_n)_{n\in\bz}$ with the following properties:
\begin{enumerate}
\item $V_n\subset V_{n+1},\quad(n\in\bz);$ \item
$UV_{n+1}=V_n,\quad(n\in\bz);$ \item $\cup V_n$ is dense in $H$;
\item $\cap V_n=\{0\}$; \item There exists $\varphi$ in $V_0$ such
that $\{T^k\varphi\,|\,k\in\bz\}$ is an orthonormal basis for
$V_0$.
\end{enumerate}
Such a vector $\varphi$ is called an {\em orthogonal scaling
vector}. For the details of the multiresolution construction we
refer to \cite{Dau92} and \cite{BDP}. We recall here some facts.
\par
If a multiresolution is given, the wavelets can be constructed by
looking for $N-1$ functions $\psi_1,...,\psi_{N-1}$ such that
$$\{T^k\psi_i\,|\,k\in\bz,i\in\{1,...,N-1\}\}$$
is an orthonormal basis for $V_1\ominus V_0$.
\par
If the spectral measure of the operator $T=\pi(t)$ is absolutely
continuous with respect to the Lebesgue measure on
$\bt:=\{z\in\bc\,|\,|z|=1\}$, then, by Borel functional calculus,
we can define a representation of $\linft$ on $H$ by
$$\pi(f)=f(T),\quad (f\in\linft),$$
and this representation will satisfy
$$U\pi(f)U^{-1}=\pi(f(z^N)),\quad(f\in\linft),\quad\pi(z)=T.$$
(On $\bt$, we have the Lebesgue measure $\mu$.)
\par
Since $U\varphi\in V_{-1}\subset V_0$ the scaling vector $\varphi$
satisfies the following {\em scaling equation}:
\begin{equation}\label{eq0_1}U\varphi=\sum_{k\in\bz}a_kT^k\varphi,
\end{equation}
for some coefficients $a_k\in\bc$. This can be rewritten as
\begin{equation}\label{eq0_2}
U\varphi=\pi(m_0)\varphi, \end{equation} with
$m_0(z)=\sum_{k\in\bz}a_kz^k,$ $z\in\bt$. $m_0$ is called the {\em
low-pass filter}.
\par
A vector $\varphi$ that satisfies (\ref{eq0_2}) for some
$m_0\in\linft$ and which is cyclic for the representation of the
group is called a {\em scaling vector} with filter $m_0$.
\par
Each pair of vectors $v_1,v_2$ has a {\em correlation function}
$h_{v_1,v_2}\in\lonet$ associated to it; this is defined by
$$\ip{T^kv_1}{v_2}=\int_{\mathbb{T}}z^kh_{v_1,v_2}\,d\mu,\quad(k\in\bz),$$
or, equivalently, by the Radon-Nikodym derivative of the functional $$f\mapsto\ip{\pi(f)v_1}{v_2}.$$
We denote by $h_v:=h_{v,v}$.
\par
If $\varphi$ is a scaling vector with filter $m_0$, then its
correlation function satisfies the following equation
$$R_{m_0}h_\varphi=h_\varphi,$$
where $R_{m_0}$ is the Ruelle operator defined on $\lonet$ by
$$R_{m_0}f(z)=\frac{1}{N}\sum_{w^N=z}|m_0(w)|^2f(w),\quad(f\in\lonet,z\in\bt).$$
\par
The interesting fact is that the converse is also true, in the
sense that each scaling equation (\ref{eq0_2}) has a solution in
some Hilbert space such that the correlation function of the
scaling vector is some prescribed function $h$ with $R_{m_0}h=h$.
More precisely, we have:
\begin{theorem}\label{th0_1}\cite{Jor01},\cite{Dut3}.
Let $m_0\in\linft$ and $h\in\lonet$ such that $h\geq 0$ and
$R_{m_0}h=h$. Then there exist a representation $\pi$ of the group
$BS(1,N)$ on a Hilbert space $H$ and $\varphi\in H$ such that
$$U\varphi=\pi(m_0)\varphi,$$
and the correlation function of $\varphi$ is $h$. Moreover, this
is unique up to isomorphism.
\end{theorem}
We call this representation, the {\em wavelet representation}
associated to $m_0$ and $h$. \par If $h=1$ then the scaling vector
is an orthogonal one and $m_0$ satisfies the relation
$$R_{m_0}1=1$$
which is known in the literature as the {\em quadrature mirror
filter (QMF)} equation.
\par
Our goal is to analyze these representation, perform the direct
integral decomposition and evaluate the consequences that these
have on the scaling equation.
\par
In section \ref{irreducible} we construct the theoretical
framework for this purpose. We define a special type of
representation (definition \ref{def2_1}), and we prove in section
\ref{wav} that the wavelet representations are a particular case
(theorem \ref{th3_3}). We classify these representation
(proposition \ref{prop2_3}) and give a criterion that identifies
them among the representations of the Baumslag Solitar group
(proposition \ref{prop2_3_1}).
\par
In section \ref{wav} we establish a connection between wavelet
representation and the measures associated to some random walks on
$\mathbb{T}$ which were introduced in a recent paper by Palle
Jorgensen \cite{Jor04}. Theorem \ref{th3_3} is the central result
of the paper and it shows that the wavelet representations can be
realized on the solenoid $\SN$ using the measure associated to a
random walk. Corollary \ref{cor3_6} establishes the 1-1
correspondence between operators in the commutant, functions which
are invariant for the shift $S$, and fixed points of the transfer
operator $R_{m_0}$.
\par
In theorem \ref{th3_6} we give the direct integral decomposition
of the wavelet representation which corresponds to the
decomposition of the random walk measure into ergodic components.
\par
Section \ref{ex} is dedicated to examples. In \ref{ltwo} we
analyze the representation on $\ltwor\oplus...\oplus\ltwor$ (which
obviously contains the classical case on $\ltwor$). We describe
the direct integral decomposition (theorem \ref{th4_2_1}). This
generalizes the result from \cite{LPT}. Also we describe the
associated measure $m$ on the solenoid extending in this way some
results from \cite{Jor04}.
\par
As some consequences of the decomposition we mention theorem
\ref{th4_4} and corollary \ref{cor4_5} which give general
necessary and sufficient conditions for the existence of
$\ltwor$-solutions of the scaling equation and also can be used as
some ergodic statements about QMF filters.
\par
In section \ref{frac} we deal with the representation on a fractal
measure that was introduced in \cite{DutJo}. We compute the
Fourier coefficients of the associated measure (proposition
\ref{prop4_2_1}), which gives also a geometric insight into the
Cantor set (lemma \ref{lem4_2_2}), and we propose a Fourier
transform for this fractal measure (corollary \ref{cor4_2_3}).
\par
The last example is for $m_0=1$ when we reobtain the Haar measure
on $\SN$ and its ergodic properties relative to the shift which
imply the irreducibility of the representation.

\section{\label{rep} Representations that have wavelets}
In this section we analyze some restrictions that are imposed on a
representation of the Baumslag Solitar group $BS(1,N)$ in the case
when this representation has a wavelet.
\par
It is known that the existence of a wavelet establishes some strong
restrictions on the representation. Some of these restrictions are
analyzed in \cite{Web} for each of the operators $U$ and $T$ in
part, for example, it is shown that $U$ and $T$ must be both
bilateral shifts of infinite multiplicity. We analyze here $U$ and
$T$ coupled by the commuting relation $UTU^{-1}=T^N$.
\par
The next results generalizes theorem 5.1 in \cite{MV00}, but the
proof follows the same ideas.
\begin{theorem}\label{th1_1}
Let $\pi$ be a representation of the group $BS(1,N)$ that has a
wavelet. Then $\pi$ extends to a faithful representation of
$C^*_r(BS(1,N))$ -the reduced $C^*$-algebra of $BS(1,N)$, i.e.,
the $C^*$-algebra generated by the right regular representation.
Each element in the $C^*$-algebra generated by $\pi$ has a
connected spectrum.
\end{theorem}
\begin{proof} Let $\{\psi_1,...,\psi_p\}$ be an orthonormal wavelet.
Since the group is amenable, $\pi$ is weakly contained in the
right regular representation. For the converse, take $g\in
BS(1,N)$. Then $g$ is of the form $g=u^jt_{\lambda}$ for some
$j\in\bz$ and some $\lambda\in\bz[1/N]$. Then, for $k\in\bz$ we
have
$$\ip{\pi(g)\pi(u^{-k})\psi_1}{\pi(u^{-k})\psi_1}=\ip{\pi(u^ju^kt_{\lambda}u^{-k})\psi_1}{\psi_1}=\ip{\pi(u^jt_{\lambda
N^k})\psi_1}{\psi_1}.$$ But, for $k$ big enough, $\lambda N^k$ is
an integer $l$, and therefore
$$\ip{\pi(g)\pi(u^{-k})\psi_1}{\pi(u^{-k})\psi_1}=\ip{\pi(u^jt^l)\psi_1}{\psi_1}=\delta_{e,g}.$$
This implies that the right regular representation is weakly
contained in $\pi$ and the first assertion is proved.
\par
Since the representation is faithful and the spectrum is invariant
under isomorphisms, the last assertion follows from lemma 2.1 in
\cite{MV00}.
\end{proof}

\section{\label{irreducible}A class of irreducible representations}
For a measure $\nu$ on $\SN$, we denote by $\nu\circ S$ the
measure defined by $\nu\circ S(E)=\nu(S(E))$ for all measurable
subsets $E$ of $\SN$. Alternatively, for $f$ measurable and
positive on $\SN$,
$$\int_{\SN}f\,d\nu\circ S=\int_{\SN}f\circ
S^{-1}\,d\nu.$$
\begin{definition}\label{def2_1}
Let $\nu$ be a $\sigma$-finite Borel measure on $\SN$ such that
$\nu\circ S$ and $\nu$ are mutually absolutely continuous; we say
that $\nu$ is quasi-invariant for $S$. Also, consider a measurable
map $\theta:\SN\rightarrow\mathbb{T}$. Let
$$\Delta:=\frac{d(\nu\circ S)}{d\nu}.$$
\par
Define $H=L^2(\nu)$,
$$U\xi=\sqrt{\Delta}\theta\xi\circ S,\quad(\xi\in
L^2(\nu)),$$
$$T\xi(z_n)_{n\in\Zu}=z_0\xi(z_n)_{n\in\Zu},\quad(\xi\in
L^2(\nu),(z_n)_{n\in\Zu}\in\SN).$$
\end{definition}
\begin{proposition}\label{prop2_2}
\begin{enumerate}
\item
 With $\nu$ and $\theta$ as in definition
\ref{def2_1}, $U$ and $T$ define a representation $\pi_{\nu,\theta}$ of the group
$BS(1,N)$. \item If $f$ is in $L^{\infty}(\nu)$ and
$n\in\mathbb{Z}$ then
$$U^nM_fU^{-n}=M_{f\circ
S^n},$$ where, for $g\in L^{\infty}(\nu)$,  $M_g$ is the
multiplication operator by $g$, $M_g\xi=g\xi$. For
$n,k\in\mathbb{Z}$, $U^{-n}T^kU^{n}$ is the multiplication by the
character $k/N^n$.
\par
The von Neumann algebra generated by
$\{U^nT^kU^{-n}\,|k,n\in\mathbb{Z}\}$ is $$\{M_f\,|\,f\in
L^{\infty}(\nu)\}.$$ \item The commutant of the representation is
given by
$$\{M_f\,|\,f\in L^{\infty}(\nu), f\circ S=f, \nu\mbox{-a.e.}\},$$
so the representation is irreducible if and only if $S$ is ergodic
with respect to $\nu$, i.e., the sets $A$ which are invariant with
respect to $S$ have $0$ or full measure
\end{enumerate}
\end{proposition}
\begin{proof}
Since $\Delta=d(\nu\circ S)/d\nu$, it follows that
$$\int_{\SN}\xi\,d\nu=\int_{\SN}\Delta\xi\circ S\,d\nu,\quad(\xi\in
L^1(\nu)),$$ and this implies that $U$ is an isometry. Since $\nu$
is also absolutely continuous with respect to $\nu\circ S$, the
inverse of $U$ is well defined by:
$$U^{-1}\xi=\frac{1}{\theta\circ S^{-1}\sqrt{\Delta}\circ S^{-1}}\xi\circ
S^{-1},\quad(\xi\in L^2(\nu)),$$ and therefore $U$ is unitary. $T$
is just a multiplication operator by a function which has absolute
value $1$ so $T$ is also unitary.
\par
Now take $n,k\in\mathbb{Z}$ and compute by induction
\begin{equation}\label{eq2_2_1}
U^n\xi=(\theta\sqrt{\Delta})^{(n)}\xi\circ S^n,
\end{equation}
where, for a function $f$ on $\SN$ we use the notation:
\begin{equation}
\label{eq2_2_2}
f^{(n)}=\left\{\begin{array}{ccc}
f\circ S...f\circ S^{n-1},&\mbox{ if }&n>0,\\
1&\mbox{ if }&n=0,\\
\frac{1}{f\circ S^{-1}...f\circ S^{-n}}&\mbox{ if }&n<0.
\end{array}\right.
\end{equation}
Then, after a straightforward computation, we obtain:
$$U^{-n}T^kU^{n}\xi(z_i)_{i\in\mathbb{Z}}=z_n^k\xi(z_i)_{i\in\mathbb{Z}},$$
and if $f$ is in $L^{\infty}(\SN)$ then
$$U^{n}M_fU^{-n}=M_{f\circ S^n}.$$
In particular $UTU^{-1}=T^N$, so this is indeed a representation
of $BS(1,N)$.
\par
The proof of the last statement in (ii) is a standard argument of
the duality theory: by Pontrjagin's duality theory, the linear
span of characters $k/N^n$ is uniformly dense in $C(\SN)$ and any
function in $L^{\infty}(\nu)$ can be approximated pointwise $\nu$
a.e. by continuous ones, so the von Neumann algebra generated by
the operators of multiplication by characters is $L^{\infty}(\SN)$
(we will use the identification between a function
$f\in L^{\infty}(\SN)$ and the multiplication operator $M_f$).
\par
It remains to compute the commutant. If $S$ is an operator that
commutes with $U$ and $T$, it must commute with the entire von
Neumann algebra generated by the elements of the form
$U^nT^kU^{-n}$, and we saw that this is $L^{\infty}(\nu)$. Since
this algebra is maximal abelian, $S$ must belong to it, so $S=M_f$
for some $f\in L^{\infty}(\nu)$. However, $S$ must commute with
$U$, too, so $M_f=UM_fU^{-1}=M_{f\circ S}$. Therefore $f=f\circ
S$. $S$ is ergodic if and only if the only such functions are the
ones that are constant $\nu$-a.e. In conclusion, the commutant is
trivial and the representation is irreducible if and only if $S$
is ergodic.
\end{proof}
\begin{proposition}\label{prop2_3}
Let $(\nu_1,\theta_1)$ and $(\nu_2,\theta_2)$ be as in the
definition \ref{def2_1}, and let $U_1,T_1$ and $U_2,T_2$ the
corresponding representations. The representations are equivalent
if and only if $\nu_1$ and $\nu_2$ are mutually absolutely
continuous and $\theta_1$ and $\theta_2$ are cocycle equivalent,
in the sense that there exists a function
$\lambda:\SN\rightarrow\mathbb{T}$ such that
$$\lambda\theta_1=\lambda\circ S\theta_2,\nu_1\mbox{-a.e.}$$
\end{proposition}
\begin{proof}
Suppose the representations are equivalent and let
$W:L^{2}(\nu_1)\rightarrow L^{2}(\nu_2)$ be an intertwining
isomorphism. Then, restricting our attention to $\bz[1/N]$, we see
that $W$ establishes the equivalence between the two
representations of the abelian algebra $C(\SN)$ by multiplication
operators on $\nu_1$ and $\nu_2$. This implies that $\nu_1$ and
$\nu_2$ are mutually absolutely continuous.
\par
Let $\eta:=d\nu_1/d\nu_2$. The operator $\tilde
W:L^{2}(\nu_1)\rightarrow L^{2}(\nu_2)$ defined by
$$\tilde W\xi_1=\sqrt{\eta}{\xi_1},\quad(\xi_1\in
L^{\infty}(\nu_1)).$$ is an isomorphism. Then, $\tilde W W^*$ is a
unitary operator which commutes with $L^{\infty}(\nu_2)$.
Therefore, as this is a maximal abelian subalgebra, $\tilde
WW^*=M_{\lambda}$ for some function $\lambda$ on $\SN$ which has
absolute value $1$ a.e. Then $W=M_{\cj\lambda}M_{\sqrt{\eta}}$.
\par
Since $WU_1=U_2W$, we get for $\xi\in L^2(\nu_1)$,
$$\cj{\lambda}\sqrt{\eta}\theta_1\sqrt{\Delta_1}\xi_1\circ
S=\theta_2\sqrt{\Delta_2}\cj{\lambda}\circ S\sqrt{\eta\circ
S}\xi_1\circ S.$$ But
\begin{equation}\label{eq2_3_1}
\eta\circ S=\frac{d\nu_1\circ S}{d\nu_2\circ S}= \frac{d\nu_1\circ
S}{d\nu_1}\frac{d\nu_1}{d\nu_2}\frac{d\nu_2}{d\nu_2\circ
S}=\Delta_1\eta\frac{1}{\Delta_2},
\end{equation}
so $\cj{\lambda}\theta_1=\cj{\lambda}\circ S\theta_2$,
$\nu_1$-almost everywhere.
\par
For the converse, take $\eta=d\nu_1/d\nu_2$ and define $W$ from
$L^2(\nu_1)$ to $L^2(\nu_2)$ by
$$W\xi_1=\lambda\sqrt{\eta}\xi_1,\quad(\xi_1\in L^2(\nu_1)).$$
This defines an isomorphism which clearly intertwines $T_1$ and
$T_2$, and using equation (\ref{eq2_3_1}), we see that it also
intertwines $U_1$ and $U_2$.
\end{proof}
\begin{proposition}\label{prop2_3_1}
Let $\pi$ be a representation of the group $BS(1,N)$ such that the
restriction of the representation to the subgroup $\Zu[1/N]$ has a
cyclic vector. Then there is a quasi-invariant probability measure
$\nu$ on $\SN$, and a map $\theta:\SN\rightarrow\mathbb{T}$ such
that $\pi$ is equivalent to the representation $\pi_{\nu,\theta}$ (see Proposition \ref{prop2_2}).
\end{proposition}
\begin{proof}
By the Stone-Mackey theorem applied to the abelian group
$\Zu[1/N]$ we can find a measure $\nu$ on the dual group $\SN$ and
a measurable multiplicity function
$m:\SN\rightarrow\{0,1,...,\infty\}$ such that there is an isometric isomorphism $\Phi$ from the Hilbert space of the representation $H$ to 
$L^2(\SN,\nu,m):=\oplus_{j\geq1}L^2(\{z\in\SN\,|\,m(z)\geq j\},\nu)$ which transforms the representation into multiplication operators, 
$$(\Phi\pi(\lambda)\Phi^{-1})(\xi_j)_{j\geq1}=(\chi_{\lambda}\xi_j)_{j\geq1},\quad((\xi_j)_j\in L^2(\SN,\nu,m)).$$
Here $\chi_{\lambda}$ is the character of the group $\SN$ given by the duality in (\ref{eqdual}):
\begin{equation}\label{eqchar}
\chi_{\lambda}((z_n)_{n\in\bz})=\ip{\lambda}{(z_n)_{n\in\bz}},\quad((z_n)_{n\in\bz}\in\SN).\end{equation}
\par
Since there is a cyclic vector, the
multiplicity function can be taken to be constant $1$, and we can
take $\nu(\SN)=1$. Thus the representation of $\Zu[1/N]$ is
equivalent to the representation on $L^2(\nu)$ by multiplications
by characters, and $\Phi\pi(\lambda)\Phi^{-1}=M_{\chi_{\lambda}}$, for
$\lambda\in\Zu[1/N]$. So (by composition with the isomorphism $\Phi$) we can assume that the representation $\pi$ is on $L^2(\nu)$ and $\pi(\lambda)=M_{\chi_\lambda}$.
\par
The representation $\pi(u)$ of the other generator $u$ of the group $BS(1,N)$ is a unitary $U$ with the property
$$U\pi(\lambda)U^{-1}=\pi(u\lambda u^{-1})=\pi(N\lambda),\quad(\lambda\in\Zu[1/N]).$$
But this implies, by approximation, that
$$UM_fU^{-1}=M_{f\circ S},\quad(f\in L^{\infty}(\nu)).$$
Then $U\xi=\xi\circ S U1$, for $\xi\in L^{\infty}(\nu)$. Denote by
$f:=U1$ then, since $U$ is unitary, we have that, for $\xi\in
L^{\infty}(\nu)$,
$$\int_{\SN}|\xi\circ S^{-1}|^2\,d\nu=\int_{\SN}|f|^2|\xi|^2\,d\nu,$$
so $d(\nu\circ S)/d\nu=|f|^2$, and $f$ does not vanish on a set of
positive measure, so $\nu$ is quasi-invariant. Take $\theta=f/|f|$
and everything follows.
\end{proof}

\begin{example}\label{ex2_4}
Take a point $x=(z_i)_{i\in\Zu}\in\SN$ and let $\nu$ be the
counting measure on the orbit on this point under $S$:
$$\nu(E)=\mbox{card}(\{S^n(x)\,|\,n\in\mathbb{Z}\}\cap E),\quad
(E\subset\SN).$$ Take $\theta=1$.
\par
$\nu$ and $\theta$ satisfy the requirements of proposition
\ref{prop2_2} and let $U$, $T$ be the corresponding
representation. This representation is irreducible because $S$ is
ergodic with respect to the measure $\nu$. We distinguish two
cases
\begin{enumerate}
\item If $x$ is periodic (which means that the orbit is finite),
i.e., there exists $p>0$ such that $z_{i+p}=z_i$ for all
$i\in\mathbb{Z}$, then we can take $p$ minimal with this property
and we see that the representation is equivalent to the following
on $H_x=\mathbb{C}^p$:
$$U_x(\xi_0,...,\xi_{p-1})=(\xi_{p-1},\xi_0,...,\xi_{p-2}),$$
$$T_x(\xi_0,...,\xi_{p-1})=(z_0\xi_0,z_1\xi_1,...z_{p-1}\xi_{p-1}).$$
\item If $x$ is not periodic, then the representation can be
realized on $l^2(\Zu)$, with
$$U_x\xi(k)=\xi(k-1),\quad(k\in\mathbb{Z}),$$
$$T_x\xi(k)=z_k\xi(k),\quad(k\in\mathbb{Z}).$$
\end{enumerate}
\par
The representations associated to two points $x,x'\in\SN$ are
equivalent if and only if $x$ and $x'$ are on the same orbit.

\par
These representations can be obtained also using the Mackey
machine (see for example chapter 6 of \cite{Fol}) as some induced
representations. However, the technique of Mackey does not give
all the irreducible representations because the action of $\Zu$ on
$\Zu[1/N]$ is not regular (see theorem 6.42 in \cite{Fol}). One
example of an irreducible representation that does not come from
the Mackey construction is given next.
\end{example}

\begin{example}\label{ex2_6}
Consider $\mu_{\SN}$ the Haar measure on $\SN$. By the uniqueness
of the Haar measure, $\mu_{\SN}\circ S=\mu_{\SN}$. By a theorem of
Rohlin-Halmos (see \cite{Wal}), $S$ is ergodic with respect to the
Haar measure if and only if the only character $\lambda=k/N^p$
satisfying $\alpha_n(\lambda)=\lambda$ for some
$n\in\mathbb{Z}\setminus\{0\}$, is the trivial one. But this is
clear because $\alpha_n(\lambda)=N^n\lambda.$
\par
This implies that the representation associated to the Haar
measure as in definition \ref{def2_1} is irreducible. We will see
that this is actually a wavelet representation associated to the
filter $m_0=1$.
\end{example}

\section{\label{wav}Wavelet representations and random walks}
\par
In a recent paper \cite{Jor04}, Palle Jorgensen, extending some
earlier work by Richard Gundy \cite{Gun00}, has realized a
connection between wavelet theory and some probability measures
associated to certain random walks. We recall here the definition of these
measures and refer the reader to \cite{Jor04} for the details. The
measures are perfectly adapted to our purpose and we prove in
theorem \ref{th3_3} that the wavelet representation associated to
some filter $m_0$ can be realized on such a measure. \par Consider
the $N$ inverse branches of the map $\sigma:x\mapsto Nx\mod 1$, on
$[0,1)$, $\tau_k:[0,1)\rightarrow[\frac{k}{N},\frac{k+1}{N})$,
$$\tau_k(x)=\frac{x+k}{N},\quad (x\in[0,1],k\in\{0,...,N-1\}).$$
\par
Denote by $\Omega:=\{0,...,N-1\}^{\mathbb{N}}$ with the product
topology.
\par
Consider a function $W\in L^\infty[0,1]$, $W\geq 0$, with the
property that
\begin{equation}\label{eq3_0_1}
\sum_{k=0}^{N-1}W(\tau_k(x))=1,\quad(x\in[0,1)).
\end{equation}
For example, if $m_0\in\linft$ satisfies $R_{m_0}1=1$ then $W$
defined by $W(x)=|m_0(e^{-2\pi\mbox{i}x})|^2/N$ will satisfy
(\ref{eq3_0_1}).
\par
We can identify functions on $[0,1)$ with functions on $\bt$ by
$$W(x)\leftrightarrow W(e^{-2\pi\mbox{i}x}).$$
Also, we can identify functions $f$ on $\mathbb{T}$ with functions
on the $\SN$ that depend only on the first coordinate:
$$f((z_n)_{n\in\bz}):=f(z_0),\quad((z_n)_{n\in\bz}\in\SN).$$
\par
It is proved in \cite{Jor04}, that for each $x\in[0,1)$ there
exists a probability measure $P_x$ on $\Omega$ such that, if a
function $f$ on $\Omega$ depends only on a finite number of
coordinates, $\omega_1,...,\omega_n$, then
\begin{equation}\label{eq3_1}
\int_{\Omega}fdP_x=\sum_{\omega_1,...,\omega_n}f(\omega_1,...,\omega_n)W(\tau_{\omega_1}x)W(\tau_{\omega_2}\tau_{\omega_1}x)...W(\tau_{\omega_n}...\tau_{\omega_1}x).
\end{equation}
\par
We can identify the space $[0,1)\times\Omega$ with the solenoid
$\SN$:
\begin{proposition}\label{prop3_1}
The map $\Phi:[0,1)\times\Omega\rightarrow\SN$ defined by
$$\Phi(x,\omega)=(e^{-2\pi\mbox{i}x},e^{-2\pi\mbox{i}\tau_1x},e^{-2\pi\mbox{i}\tau_2\tau_1x},...,e^{-2\pi\mbox{i}\tau_n...\tau_1x},...),\quad(x\in[0,1),\omega\in\Omega),$$
is a measurable bijection.
\end{proposition}
Note that under this identification, the shift takes the form
$$\Phi^{-1}S\Phi: (x,\omega)\mapsto
(\sigma(x),i_x\omega_1\omega_2...),\mbox{ where }i_x=k\mbox{ if
}x\in[\frac{k}{N},\frac{k+1}{N}),$$ and its inverse is
$$\Phi^{-1}S^{-1}\Phi:
(x,\omega)\mapsto(\tau_{\omega_1}(x),\omega_2\omega_3,...)$$
\par
Define the measure $m$ on $\SN$ by
\begin{equation}\label{eq3_2}
\int_{\SN}f\,dm=\int_{[0,1)}\int_{\Omega}f(\Phi(x,\omega))\,dP_x(\omega)\,
dx,\quad(f\in C(\SN)).
\end{equation}
\begin{proposition}\label{prop3_2}
$m$ is a probability measure on $\SN$ with the following
properties:
\begin{enumerate}
\item  If $f\in L^1(m)$ depends only on the first $n$ coordinates
$z_0,...,z_{n-1}$ then
$$\int_{\SN}f\,dm=\int_{\mathbb{T}}\sum_{w^{N^n}=z}f(w^{N^{n-1}},w^{N^{n-2}},...,w^N,w)W^{(n)}(w)\,dz,$$
where $W^{(n)}(z)=W(z)W(z^N)...W(z^{N^{n-1}}).$ \item $m$ is the
unique probability measure on $\SN$ which satisfies the conditions
\begin{equation}\label{eq3_2_1}
\int_{\SN}f\,dm=\int_{\mathbb{T}}f(z)\,dz,\quad(f\in
L^1(\mathbb{T})), \end{equation} and
\begin{equation}\label{eq3_2_2}
\int_{\SN}f\circ S^{-1}\,dm=\int_{\SN}NWf\,dm,\quad(f\in L^1(m)),
\end{equation}
(i.e. $d(m\circ S)/dm=NW$.) \item For $f\in L^1(m)$, and
$n\geq0$, (also for $n<0$ when $W$ does not vanish on a set of positive measure):
\begin{equation}\label{eq3_2_3}
\int_{\SN}f\circ S^{-n}\,dm=\int_{\SN}N^nW^{(n)}f\,dm.
\end{equation}
\item If
$W(z)=\sum_{k\in\Zu}a_kz^k$ then 
$$\hat
m(\lambda):=\int_{\SN}\chi_{\lambda}\,dm=\int_{\mathbb{T}}z^lN^pW^{(p)}\,dz,\quad(\lambda=l/N^p\in\Zu[1/N]),$$
(where $\chi_{\lambda}$ is the character on $\SN$ attached to
$\lambda$, see (\ref{eqchar})). Moreover $\hat m(k)=\delta_k$ for $k\in\Zu$ and $\hat m$
satisfies the following scaling equation:
\begin{equation}\label{eq3_2_4}
\hat m(\lambda)=N\sum_{k\in\Zu}a_k\hat
m(N\lambda+k),\quad(\lambda\in\Zu[1/N]).
\end{equation}
\end{enumerate}
\end{proposition}
\begin{proof}
If $f\in L^1(m)$ depends only on the first $n$ coordinates, then
$f\circ\Phi$ depends only on $x$ and $\omega_1,...,\omega_{n-1}$.
Also,
$$f\circ\Phi(x,\omega_1,...,\omega_{n-1})=f(w^{N^{n-1}},...,w),\mbox{ with
}w=e^{-2\pi\mbox{i}\tau_{\omega_{n-1}}...\tau_{\omega_1}(x)}.$$ So (i)
follows from (\ref{eq3_1}).
\par
Equation (\ref{eq3_2_1}) is clear and, to prove (\ref{eq3_2_2}),
take $f\in C(\SN)$ which depends only on the first $n$
coordinates. Then $f\circ S^{-1}$ depends only on the first $n+1$
coordinates and
\begin{align*}
\int_{\SN}f\circ
S^{-1}\,dm&=\int_0^1\sum_{\omega_1,...,\omega_n}f\circ
S^{-1}\circ\Phi(x,\omega_1,...,\omega_n)W(\tau_{\omega_1}x)...W(\tau_{\omega_n}...(\tau_{\omega_1}x))\,dx\\
&=\int_0^1\sum_{\omega_1,...,\omega_n}f\circ\Phi(\tau_{\omega_1}x,\omega_2,...,\omega_n)W(\tau_{\omega_1}x)...W(\tau_{\omega_n}...(\tau_{\omega_1}x))\,dx\\
\end{align*}
If we denote by
$$g(x):=\sum_{\omega_2,...,\omega_n}f\circ\Phi(x,\omega_2,...,\omega_n)W(\tau_{\omega_2}x)...W(\tau_{\omega_n}...(\tau_{\omega_2}x))$$
and use the fact that (with a change of variable)
$$\int_0^1\sum_{\omega_1}W(\tau_{\omega_1}x)g(\tau_{\omega_1}x)\,dx=\int_0^1NW(x)g(x)\,dx,$$
then we obtain
$$\int_{\SN}f\circ S^{-1}\,dm=\int_0^1NW(x)\sum_{\omega_2,...,\omega_n}f\circ\Phi(x,\omega_2,...,\omega_n)W(\tau_{\omega_2}x)...W(\tau_{\omega_n}...(\tau_{\omega_2}x))\,dx$$
$$=\int_{\SN}NWf\,dm.$$
\par
(iii) follows from (ii) by induction. All functions on $\SN$ can be approximated by functions of the form $f\circ S^{-n}$ with $f\in L^1(\bt)$ and $n\geq0$ (these are the functions which depend only on the first $n$ coordinates). 
 If a measure $m'$ satisfies the conditions of (ii), then
$$\int_{\SN}f\circ S^{-n}\,dm'=\int_{\SN}N^{n}W^{(n)}f\,dm'=$$
$$\int_{\mathbb{T}}N^{n}W^{(n)}(z)f(z)\,dz.$$
Thus the conditions of (ii) determine $m$ uniquely.
\par
For (iv), observe that $\chi_{\lambda}\circ
S^{p}=\chi_{N^p\lambda}$, for any $\lambda\in\Zu[1/N]$ and any
$p\in\Zu$. Then, using (iii),
$$\hat
m(l/N^p)=\int_{\SN}N^pW^{(p)}\chi_{l}\,dm=\int_{\mathbb{T}}z^lN^pW^{p}.$$
\par
From (\ref{eq3_2_2}), with $f=\chi_{\lambda}$, we obtain
$$\hat
m(\lambda/N)=\int_{\mathbb{T}}N\sum_{k\in\Zu}a_kz_0^k\chi_\lambda\,dm(z_n)_{n\geq0}=N\sum_{k\in\Zu}a_k\int_{\mathbb{T}}\chi_{\lambda+k}\,dm$$
which implies (\ref{eq3_2_4}).
\end{proof}

\begin{theorem}\label{th3_3}
Let $m_0\in\linft$ be non-singular (i.e., it does not vanish on
a set of positive measure), with $R_{m_0}1=1$. Let
$W:=|m_0|^2/N$ and $\theta=m_0/|m_0|$. Let $m$ be the measure
associated to $W$. Then the wavelet representation associated to
$m_0$ and $h=1$ is the representation of the Baumslag Solitar
group $BS(1,N)$ associated to $m$ and $\theta$, i.e., $H=L^2(m)$,
$$U\xi=m_0\xi\circ S,\quad \pi(f)\xi=f\xi,\quad(\xi\in
L^2(m),f\in\linft),$$
$$\varphi=1.$$
The commutant of the representation is
$$\{M_f\,|\, f\in L^{\infty}(m),f\circ S=f\, \nu\mbox{-a.e.}\}.$$
\end{theorem}
\begin{proof}
By proposition \ref{prop3_2} (ii), we know that $d(m\circ
S)/dm=NW=|m_0|^2$, and since $m_0$ is non-singular, we have also
that $m$ is absolutely continuous with respect to $m\circ S$.
Therefore we can use proposition \ref{prop2_2} and we obtain the
representation of $BS(1,N)$. With (\ref{eq3_2_1}) we see that, by
Borel functional calculus, $T$ generates a representation $\pi$ of
$\linft$, $\pi(f)=f(T)=M_f$ for $f\in\linft$.
\par
Also, note that
$$\ip{\pi(f)\varphi}{\varphi}=\int_{\SN}f\,dm=\int_{\mathbb{T}}f(z)\,dz,$$
$$U\varphi=\pi(m_0)\varphi.$$
Also
$$U^{-n}\pi(f)U^n\varphi=M_{f\circ S^{-n}}\varphi=f\circ S^{-n},$$
and, because the functions that depend on only finitely
many coordinates are dense in $L^2(m)$, it follows that $\varphi$
is cyclic for the representation, so it is a scaling vector.
\par
The commutant is obtained from proposition \ref{prop2_3}.
\end{proof}

\begin{proposition}\label{prop3_4}
The map $E$ from $L^1(m)$ to $L^1(\mathbb{T})$ defined by
$$E(f)(e^{-2\pi
ix})=\int_{\Omega}f\circ\Phi(x,\omega)\,dP_x(\omega),\quad(f\in
L^1(m),x\in [0,1)),$$ is a well defined conditional expectation
(i.e., $E^2=E$, $E(f)\geq 0$ if $f\geq 0$, $E(gf)=gE(f)$ if
$g\in\linft$ and $f\in L^1(m)$). Moreover,
$$\int_{\SN}f\,dm=\int_{\mathbb{T}} E(f)\,dz,\quad(f\in L^1(m)).$$
$E$ maps $L^{\infty}(m)$ into $\linft$ and
$\|E(f)\|_{\infty}\leq\|f\|_{\infty}$.
\par
$E$ maps $L^2(m)$ into $L^2(\mathbb{T})$ and the restriction of
$E$ to $L^2(m)$ coincides with the projection onto the subspace of
functions that depend only $z$ (which can be identified with
$L^2(\mathbb{T})).$
\end{proposition}
\begin{proof}
Everything can be checked by some straightforward computations.
\end{proof}
\begin{proposition}\label{prop3_5}
If $\xi_1,\xi_2\in L^2(m)$, then their correlation function is
$h_{\xi_1,\xi_2}=E(\xi_1\cj{\xi}_2)$.
\end{proposition}
\begin{proof}
For $f\in\linft$ we have
$$\ip{\pi{f}(\xi_1)}{\xi_2}=\int_{\SN}f\xi_1\cj{\xi}_2\,dm=\int_{\mathbb{T}}E(f\xi_1\cj{\xi}_2)\,dz=\int_{\mathbb{T}}fE(\xi_1\cj{\xi}_2)\,dz.$$
\end{proof}
\begin{corollary}\label{cor3_6}
There is a one-to-one linear and monotone correspondence between
the following data:
\begin{enumerate}
\item Operators $S$ in the commutant of $\{U,T\}$; \item Cocycles,
i.e., functions $f\in L^{\infty}(m)$ such that $f\circ S=f$; \item
Functions $h\in\linft$ which are harmonic with respect to the
Ruelle operator, i.e., $R_{m_0}h=h$.
\end{enumerate}
From (i) to (ii) the correspondence is given in theorem
\ref{th3_3}. From (ii) to (iii), the correspondence is $f\mapsto
h=E(f)$. From (iii) to (i) the correspondence is given by theorem
3.18 in \cite{Dut3}.
\end{corollary}
\begin{proof}
Everything follows from theorem \ref{th3_3}, proposition
\ref{prop3_5}, theorem 3.18 in \cite{Dut3}, and see also theorem
2.7.1 in \cite{Jor04}.
\end{proof}

\begin{theorem}\label{th3_6}
Let $m_0\in\linft$ be a non-singular filter with $R_{m_0}1=1$, and
let $(H,U,\pi,\varphi)$ be the wavelet representation associated
to $m_0$. Then there is a standard measure space
$(A,\mathcal{M},\mu)$, a measurable field $\{\mathcal{H}_a\}$ of
Hilbert spaces on $A$, a measurable field $\{\pi_a\}$ of
irreducible representations $BS(1,N)$ and a unitary map
$\Psi:H\rightarrow\int^{\oplus}\mathcal{H}_a\,d\mu(a)$, such that
\begin{enumerate}
\item $\Psi\pi(x) \Psi^{-1}=\int^{\oplus}\pi_{a}(x),$ for $x\in
BS(1,N)$; \item $\Psi\pi'\Psi^{-1}$ is the algebra of diagonal
operators on $\int^{\oplus}\mathcal{H}_{a}\,d\mu({a}),$ $\pi'$
being the commutant of the representation $\pi$.
\end{enumerate}
For almost every $a\in A$ there exists a unique ergodic,
quasi-invariant probability measure $\nu_a$ on $\SN$ such that
$d(\nu_a\circ S)/d\nu_a=|m_0|^2$ and if $\theta=m_0/|m_0|$, then
$\pi_a$ is equivalent to the representation $\pi_{\nu_a,\theta}$.
Moreover, for $f\in C(\SN)$, if $m$ is the measure associated to
$m_0$ as in theorem \ref{th3_3}, then
$$\int_{\SN}f\,dm=\int_{A}\int_{\SN}f\,d\nu_a\,d\mu(a).$$
\end{theorem}
\begin{proof}
Since the commutant of the representation, $\pi'$, is abelian
(theorem \ref{th3_3}) the first statements follow from the well
known direct integral decomposition theory for locally compact
groups (see theorem 7.37 and 7.38 in \cite{Fol}).
\par
We want to find now $\nu_a$. We prove first that for almost every
$x$ the vector $\Psi\varphi(a)$ is cyclic for the restriction of
the representation $\pi_a$ to the subgroup $\Zu[1/N]$. Suppose
not. Then we can find a subset $E$ of positive measure such that,
for $a$ in $E$,
$$K_a:=\overline{\mbox{span}}\{\pi_a(\lambda)\Psi\varphi(a)\,|\,\lambda\in\Zu[1/N]\}\neq\mathcal{H}_a.$$
But then we can define a measurable section
$\xi:A\rightarrow\int^{\oplus}\mathcal{H}_a\,d\mu(a)$ such that
$\xi(a)=0$ for $a\in A\setminus E$ and, for $a\in E$, $\xi(a)$ has
norm $1$ and is orthogonal to $K_a$. We have that
$\xi$ is orthogonal to
$$\overline{\mbox{span}}\{\pi(\lambda)\Psi\varphi\,|\,\lambda\in\Zu[1/N]\},$$
which contradicts the fact that $\varphi$ is cyclic for
$\pi(\Zu[1/N])$ (see theorem \ref{th3_3}).
\par
Since $\Psi\varphi(a)$ is cyclic for the representation
$\pi_a(\Zu[1/N])$, by proposition \ref{prop2_3_1}, there exists a
quasi-invariant probability measure $\nu_a'$ and a map $\theta_a$
such that $\pi_a$ is equivalent to $\pi_{\nu_a',\theta_a}$.
\par
Since $U\varphi=\sum_{k\in\mathbb{Z}}a_kT^k\varphi$ it follows
that
$$\pi_a(U)\Psi\varphi(a)=\sum_{k\in\mathbb{Z}}a_k\pi_a(T^k)\Psi\varphi(a)\,\mbox{ for almost every }a\in A.$$
Let $f_a:=d(\nu_a'\circ S)/d\nu_a'$. Then we obtain, identifying
$\mathcal{H}_a$ with $L^2(\nu_a')$,
that\begin{equation}\label{eq3_6_1_0}
\theta_a\sqrt{f_a}\Psi\varphi(a)\circ
S=\sum_{k\in\mathbb{Z}}a_kz^{k}\Psi\varphi(a)=m_0\Psi\varphi(a),\,\nu_a-\mbox{a.e}.\end{equation}
Since $\Psi\varphi(a)$ is cyclic for $\pi_a(\Zu[1/N])$, which are
multiplication operators, $\Psi\varphi(a)$ can not be zero on a
set of positive $\nu_a'$-measure.
\par
Consider now the measure
$d\nu_a=\frac{1}{C_a}|\Psi\varphi(a)|^2d\nu_a'$, where
$C_a=\int_{\SN}|\Psi\varphi(a)|^2\,d\nu_a$. The measures are
mutually absolutely continuous. Also
\begin{equation}\label{eq3_6_1}
\frac{d(\nu_a\circ S)}{d\nu_a}=\frac{d(\nu_a\circ
S)}{d(\nu_a'\circ S)}\frac{d(\nu_a'\circ
S)}{d\nu_a'}\frac{d\nu_a'}{d\nu_a}=|\Psi\varphi(a)\circ
S|^2f_a\frac{1}{|\Psi\varphi(a)|^2}=|m_0|^2. \end{equation} A
simple computation shows that $\nu_a(\SN)=1$. Since the
representation $\pi_a$ is irreducible, by proposition
\ref{prop2_2}, $S$ is ergodic w.r.t $\nu_a'$ and since this is
equivalent to $\nu_a$, $S$ is ergodic w.r.t. $\nu_a$.
\par
Using equation (\ref{eq3_6_1_0}), we obtain
$$\theta_a\frac{\Psi\varphi(a)\circ S}{|\Psi\varphi(a)\circ
S|}=\frac{m_0}{|m_0|}\frac{\Psi\varphi(a)}{|\Psi\varphi(a)|},$$
hence $\theta_a$ and $\theta$ are cocycle equivalent, and
therefore, by proposition \ref{prop2_3}, the representations
$\pi_{\nu_a',\theta_a}$ and $\pi_{\nu_a,\theta}$ are equivalent.
\par
To prove the uniqueness of the measure $\nu_a$, take $\rho_a$ a
measure with the same properties. Since they generate equivalent
representations, by proposition \ref{prop2_3}, they must be
mutually equivalent. Let $\eta=d\nu_a/d\rho_a$. Then
$$|m_0|^2=\frac{d(\nu_a\circ
S)}{d\nu_a}=\frac{d(\nu_a\circ S)}{d(\rho_a\circ
S)}\frac{d(\rho_a\circ
S)}{d\rho_a}\frac{d\rho_a}{d\nu_a}=\eta\circ
S|m_0|^2\frac{1}{\eta}.$$ Therefore $\eta=\eta\circ S$, $\nu_a$
a.e. But $\nu_a$ is ergodic so $\eta$ is constant, and, since both
measures are probability measures, $\eta=1$ a.e.
\par
For the last equality, it is enough to take $f$ to be
$\chi_{\lambda}$ for some $\lambda\in\Zu[1/N]$, because these are
uniformly dense in $C(\SN)$. Then
$$\int_{\SN}\chi_{\lambda}\,dm=\ip{\pi(\lambda)\varphi}{\varphi}=\int_{A}\ip{\pi_a(\lambda)\Psi\varphi(a)}{\Psi\varphi(a)}\,d\mu(a)$$$$=
\int_{A}\int_{\SN}\chi_{\lambda}|\Psi\varphi(a)|^2\,d\nu_a'\,d\mu(a)=\int_{A}\int_{\SN}\chi_{\lambda}\,d\nu_a\,d\mu(a).$$
\end{proof}
\begin{remark}\label{rem3_6}
Corollary \ref{cor3_6} shows that the projections in the commutant
of the representation correspond to sets which are invariant for
the shift. Proposition \ref{prop2_2} shows that irreducible
representations correspond to ergodic measures. Therefore, the
direct integral decomposition of the representation into
irreducible components corresponds to the direct integral
decomposition of the measure $m$ into its ergodic components (see
e.g. \cite{GS}). We can use theorem 1.1 in \cite{GS} to obtain
more information about the measures $\nu_a$:
\begin{enumerate}
\item For every Borel subset $B$ of $\SN$,
$$m(B)=\int_A\nu_a(B)\,d\mu(a).$$
\item If $a,a'\in A$ and $a\neq a'$ then the measures $\nu_a$ and
$\nu_a'$ are mutually singular.
\end{enumerate}
\end{remark}

\begin{proposition}\label{prop3_7}
With the notations of theorem \ref{th3_6}, a vector $\varphi'\in
H$ is cyclic for the representation iff $\Psi\varphi'(a)\neq 0$
for almost all $a\in A$.
\end{proposition}
\begin{proof}
A vector $\varphi'$ is cyclic for the representation iff it is
separating for the commutant. Since the commutant consists of
diagonal operators, $\Psi\varphi'$ is cyclic iff,
$\lambda_a\Psi\varphi(a)=0$ for almost all $a$ implies
$\lambda_a=0$ for almost all $a$. But this is true iff
$\Psi\varphi'(a)\neq0$ for almost all $a$.

\end{proof}
\section{\label{ex}Examples}In this section we analyze some
examples in more detail. We are interested in the measure $m$
associated to $m_0$ and in the decomposition of the wavelet
representation. We also present some interesting consequences.
\subsection{\label{ltwo}Representations on $\ltwor^p$}
\par
Let $C$ be a cycle, $C=\{z_0,...,z_{p-1}\}$, i.e., $z_1^N=z_2,
z_2^N=z_3,...,z_{p-2}^N=z_{p-1}, z_{p-1}^N=z_0$, a periodic orbit
for the map $z\mapsto z^N$, $p$ being the length of the orbit. Let
$|\alpha_0|=...=|\alpha_{p-1}|=1$. In \cite{BDP} and \cite{Dut2},
we constructed a representation of $BS(1,N)$ on $\ltwor^p$ by
$$Uf(\xi,i)=\alpha_i\frac{1}{\sqrt{N}}f(\frac{\xi}{N},(i+1)\mbox{mod}p),\quad(f\in\ltwor^p,\xi\in\br,i\in\{0,...,p-1\}),$$
$$Tf(\xi,i)=z_if(x-1,i)\quad(f\in\ltwor^p,\xi\in\br,i\in\{0,...,p-1\}).$$
We denote this representation by $\mathfrak{R}_{C,\alpha}$. Since
the representations $\mathfrak{R}_{C,\alpha}$ are isomorphic for
$C$ fixed and $\alpha$ variable (see \cite{BDP}), we will work
mostly with the case when all $\alpha_i=1$, and use the notation
$\mathfrak{R}_C$.
\par
Taking the Fourier transform the representation becomes
$$\hat Uf(\xi,i)=\sqrt{N}f(Nx,(i+1)\mbox{mod}p),\quad(f\in\ltwor^p,\xi\in\br,i\in\{0,...,p-1\}),$$
$$\hat Tf(\xi,i)=z_ie^{-2\pi\mbox{i}x}f(x,i)\quad(f\in\ltwor^p,\xi\in\br,i\in\{0,...,p-1\}).$$
\par
We will describe the direct integral decomposition of this
representation.
\par
Each $x\in\br$ can be identified with an element in
$\mathcal{S}_N$ by $\hat i(x)=(e^{-2\pi\mbox{i}N^{-k}x})_{k\in\bz}$.
Also we can identify the cycle $C$ with an element
$z_C\in\mathcal{S}_N$,
$$z_C:=(z_{(-k)\mbox{mod}p})_{k\in\bz}.$$
Then the product $x_C:=z_C\hat i(x)$ is an element of the group
$\mathcal{S}_N$ and, as in example \ref{ex2_4}, we can construct
an irreducible representation $\pi_{x_C}$ on $\ltwoz$ by:
$$\pi_{x_C}(u)\xi(k)=\xi(k-1),\quad(\xi\in\ltwoz,k\in\bz),$$
$$\pi_{x_C}(t)\xi(k)=z_{(-k)\mbox{mod}p}e^{-2\pi\mbox{i}N^{-k}x}\xi(k),\quad(\xi\in\ltwoz,k\in\bz).$$
Let $E=(-N^p,1]\cup[1,N^p).$
\begin{theorem}\label{th4_2_1}
Let $C=\{z_0,...,z_{p-1}\}$ be a cycle. The representation
$\mathfrak{R}_C$ is isometrically isomorphic to the representation
on $\ltwoez$ given by the direct integral
$\int_E^{\oplus}\pi_{z_Cx}$.
\end{theorem}
\begin{proof}
Define $\Psi$ from $\ltwor^p$ to $\ltwoez$ by
$$\Psi(f)(x,m)=\sqrt{N^{-m}}f(N^{-m}x,(-m)\mbox{mod}p),\quad(f\in\ltwor^p,x\in E,m\in\bz).$$
We claim that $\Psi$ is an isomorphism with inverse
$$\Psi^{-1}(f)(\xi,i)=\frac{1}{\sqrt{N^{-m}}}f(x,m),\quad(f\in\ltwoez,\xi\in\br,i\in\{0,...,p-1\}),$$
where $x\in E$ and $m\in\{0,...,p-1\}$ are uniquely determined by
the equations $\xi=N^{-m}x$ and $(-m)\mbox{mod}p=i$.
\par
First, we check that $\Psi$ is an isometry
\begin{align*}
\sum_{i=0}^{p-1}\int_{\mathbb{R}}|f(\xi,i)|^2\,d\xi&=\sum_{i=0}^{p-1}\sum_{l\in\mathbb{Z}}\int_{N^{pl+i}E}|f(\xi,i)|^2\,d\xi\\
&=\sum_{i=0}^{p-1}\sum_{l\in\bz}\int_E|f(N^{pl+i}x,i)|^2N^{pl+i}\,dx\\
&=\sum_{m\in\bz}\int_E|f(N^{-m}x,(-m)\mbox{mod}p)|^2N^{-m}\,dx.
\end{align*}
The fact that $\Psi^{-1}$ has the given form follows from a
one-line computation.
\par
Next we want to see that $\Psi$ intertwines the representations.
\begin{align*}
\Psi\hat U\Psi^{-1}(f)(x,m)&=\sqrt{N^{-m}}(\hat
U\Psi^{-1}(f))(N^{-m}x,(-m)\mbox{mod}p)\\&=\sqrt{N^{-m+1}}\Psi^{-1}(f)(N^{-m+1}x,(-m+1)\mbox{mod}p)\\
&=f(x,m-1).
\end{align*}
\begin{align*}
\Psi\hat
T\Psi^{-1}(f)(x,m)&=\sqrt{N^{-m}}(\hat T\Psi^{-1}(f))(N^{-m}x,(-m)\mbox{mod}p)\\
&=z_{(-m)\mbox{mod}p}e^{-2\pi\mbox{i}N^{-m}x}\sqrt{N^{-m}}\Psi^{-1}(f)(N^{-m}x,(-m)\mbox{mod}p)\\
&=z_{(-m)\mbox{mod}p}e^{-2\pi\mbox{i}N^{-m}x}f(x,m).
\end{align*}
\end{proof}
\par
Next, we will try to describe the measure $m$ associated to a
filter $m_0$ that gives scaling vectors for $\mathfrak{R}_C$.\par
Recall the following theorem from \cite{Dut2}:
\begin{theorem}\label{th4_2}
Suppose $m_0=\sum_{k\in\mathbb{Z}}a_kz^k$ is a Lipschitz function
with finitely many zeros, with $R_{m_0}1=1$ and let
$C_i=\{z_0^i,...,z^i_{p_i-1}\}$, $i\in\{1,...,n\}$ be the
$m_0$-cycles (assume there is at least one), i.e.,
$|m_0(z^i_j)|=\sqrt{N}$ for all $i,j$. Denote by $\theta^i_k$ the
argument of $z^i_{k\mbox{mod}p_i}$ ($e^{-2\pi
i\theta^i_k}=z^i_k$), and let $\alpha^i_k=m_0(z_{k\mbox{mod}
p_i})/\sqrt{N}\in\bt$. Then
\begin{equation}\label{eq4_1}
\hat\varphi^i_k(x)=\prod_{l=1}^\infty\frac{\cj\alpha^i_{k-l}m_0\left(\frac{x}{N^l}+\theta^i_{k-l}\right)}{\sqrt{N}},\quad(x\in\br,k\in\{0,...,p_i-1\},i\in\{1,...,n\}),
\end{equation}
defines an orthonormal scaling vector
$\varphi=(\varphi^i_0,...,\varphi_{p_i-1})_{i=1,n}$ in
$\oplus_{i=1}^n\ltwor^{p_i}$ for the representation
$\oplus_{i=1}^n\mathfrak{R}_{C_i,\alpha_i}$ with
$$U\varphi=\pi(m_0)\varphi,$$
$$\ip{\pi(t)^k\varphi}{\varphi}=\delta_k,\quad(k\in\bz).$$
The last equality can be rewritten as
\begin{equation}\label{eq4_2_1}
\sum_{i=1}^n\sum_{a\in\bz}\sum_{k=1}^{p_i}|\hat\varphi_k^i(x+a-\theta^i_k)|^2=1,\quad(x\in[0,1]).
\end{equation}
\end{theorem}
\par
Note that for each $x\in [0,1]$, $a\in\bz$ and
$k\in\{0,...,p_i-1\}$ a straightforward calculation based on the
fact that $N\theta^i_l\equiv\theta^i_{l+1}\mod 1$ shows that
\begin{equation}\label{eq4_2}z(x,i,a,k):=(e^{-2\pi
i(\frac{x+a-\theta^i_k}{N^l}+\theta^i_{k-l})})_{l\in\bz}=S^{-k}(z_{C_i})\hat
i(x+a-\theta^i_k)\in\SN,
\end{equation}
and $\Phi^{-1}(z(x,i,a,k))$ is of the form $(x,\omega(x,i,a,k))$.
\begin{proposition}\label{prop4_3}
If $m_0$ satisfies the hypotheses of theorem \ref{th4_2}, with the
previous assumptions and notations, for all $x\in [0,1]$, the
measure $P_x$ is atomic, concentrated on
$$\{\omega(x,i,a,k)\,|\,a\in\bz,k\in\{0,...,p_i-1\},i\in\{1,...,n\}\},$$
$$P_x(\omega(x,i,a,k))=|\hat\varphi^i_k(x+a-\theta^i_k)|^2,\quad(a\in\bz,k\in\{0,...,p_i-1\},i\in\{1,...,n\})$$
The measure $m$ is supported on
$\cup_{C}\cup_{k=0}^{p_i-1}S^{-k}z_C\hat i(\br)$.
\end{proposition}
\begin{proof}
To evaluate the measure of the set $\{\omega(x,i,a,k))\}$, we
write it as an intersection of cylinders in $\Omega$ and, using
(\ref{eq3_1}) and taking the limit, we obtain that

\begin{equation}\label{eq4_3}P_x(\omega(x,i,a,k))=\prod_{l=1}^\infty
W(z(x,i,a,k)_l)=\prod_{l=1}^\infty\frac{|m_0\left(\frac{x+a-\theta^i_k}{N^l}+\theta^i_{k-l}\right)|^2}{N}
\end{equation}
$$= |\hat\varphi^i_k(x+a-\theta^i_k)|^2.$$

But, using (\ref{eq4_2_1}), we obtain that
$$\sum_{i=1}^n\sum_{a\in\bz}\sum_{k=1}^{p_i}P_x(\omega(x,i,a,k))=1,$$
and the rest follows.
\end{proof}
\par
We can use the decomposition theorem \ref{th4_2_1} to obtain an
interesting ergodic result for quadrature mirror filters:
\begin{theorem}\label{th4_4}
Let $C=\{z_0,...,z_{p-1}\}$ be a cycle and denote by $\theta_k$
the argument of $z_k\mbox{mod}p$, i.e., $e^{-2\pi
i\theta_k}=z_{k\mbox{mod}p}$. Let $m_0\in\linft$ such that there
exists some (not necessarily orthogonal) scaling vector $\varphi$
cyclic for the representation $\mathfrak{R}_{C}$ with
$$U\varphi=\pi(m_0)\varphi.$$
Then for almost every $x\in (-N^p,1]\cup[1,N^p)$, there exists
$k_x\in\bz$ such that
\begin{equation}\label{eq4_4_0a}
\sum_{n=-\infty}^{k_x-1}\prod_{k=n-1}^{k_x}|m_0(\frac{x}{N^k}+\theta_{-k})|^2<\infty,
\end{equation}
\begin{equation}\label{eq4_4_0b}
\sum_{n=k_x+1}^\infty\frac{1}{\prod_{k=k_x+1}^{n}|m_0(\frac{x}{N^k}+\theta_{-k})|^2}<\infty.
\end{equation}
Conversely, if $m_0\in\linft$ satisfies (\ref{eq4_4_0a}) and
(\ref{eq4_4_0b}) then there exists a vector $\varphi$ which is
cyclic for the representation $\mathfrak{R}_{C}$ and satisfies
$$U\varphi=\pi(m_0)\varphi.$$
\end{theorem}
\begin{proof}
Using the decomposition given in theorem \ref{th4_2_1} we can move
the representation to $L^2(E,\ltwoz)$. We have then that
$\varphi(x)\in\ltwoz$ and
$$U_{x_C}\varphi(x)=\pi_{x_C}(m_0)\varphi(x),$$
for almost all $x\in E$. This rewrites as
\begin{equation}\label{eq4_4_1}
\varphi(x)(k-1)=m_0((x_C)_k)\varphi(x)_k,\quad(k\in\bz).
\end{equation}
\par
Since $\varphi$ is cyclic for the representation, $\varphi(x)\neq
0$ for almost all $x\in E$ (see proposition \ref{prop3_7}). So
there must be some $k_x$ such that $\varphi(x)(k_x)\neq 0$
Iterating equation (\ref{eq4_4_1}), we obtain that,
\begin{equation}\label{eq4_4_2}
\varphi(x)(n)=\left\{\begin{array}{ccc}
\varphi(x)(k_x)\prod_{k=n-1}^{k_x}|m_0(z_{(-k)\mbox{mod}p}e^{-2\pi\mbox{i}N^{-k}x})|^2,&\mbox{if}&n<k_x\\
\frac{\varphi(x)(k_x)}{\prod_{k=k_x+1}^{n}|m_0(z_{(-k)\mbox{mod}p}e^{-2\pi
iN^{-k}x})|^2},&\mbox{if}&n>k_x.\end{array}\right.
\end{equation}
Since $\varphi(x)\in\ltwoz$ for almost all $x$ and
$\varphi(x)(k_x)\neq 0$, (\ref{eq4_4_0a}) and (\ref{eq4_4_0b})
follow.
\par
For the converse, define $$\varphi(x)(k_x)=c_x,$$ $c_x$ being some
non-zero constant that we will compute later. Define
$\varphi(x)(k)$ by equation (\ref{eq4_4_2}), for $k\in\bz$.
Equation (\ref{eq4_4_0b}) implies that there are no zeros in the
denominators, so $\varphi(x)(k)$ is well defined, it satisfies the
equation (\ref{eq4_4_1}), and using the hypothesis, $\varphi(x)$
in $\ltwoz$ for almost all $x$. Take now $c_x$ such that
$\|\varphi(x)\|_2=1$. Then we get that $\varphi\in L^2(E,\ltwoz)$.
\par
Also, from (\ref{eq4_4_1}) it follows that
$U\varphi=\pi(m_0)\varphi$.
\par
The fact that $\varphi$ is cyclic for the representation follows
from proposition \ref{prop3_7}.
\end{proof}
\begin{corollary}\label{cor4_5}
Let $m_0\in\linft$ and let $M:=\|m_0\|_\infty$. Suppose the
following conditions are satisfied: the following limit exists for
almost all $x\in (-N^p,-1]\cup[1,N^p)$, and
\begin{equation}\label{eq4_5_1}
\lim_{n\rightarrow\infty}|m_0(\frac{x}{N^k}+\theta_{-k})|>1;
\end{equation}
There exists $1>\epsilon\geq 0$ such that
\begin{equation}\label{eq4_5_2}
\mu(\{z\in\mathbb{T}\,|\, |m_0(z)|^2\leq\epsilon\})>\frac{\ln
M}{\ln \frac{M}{\epsilon}}.
\end{equation}
(When $\epsilon=0$ this means that $m_0$ is zero on a set of
positive measure.)
\par
Then there exists a scaling vector $\varphi\in\ltwor^p$ which is
cyclic for the representation $\mathfrak{R}_C$ and satisfies the
scaling equation
$$U\varphi=\pi(m_0)\varphi.$$
\end{corollary}
\begin{proof}
\end{proof}
We use theorem \ref{th4_4}. Using the ratio test, equation
(\ref{eq4_5_1}) implies equation (\ref{eq4_4_0b}).
\par
Let
$$A_\epsilon:=\{z\in\mathbb{T}\,|\, |m_0(z)|^2\leq\epsilon\},\quad\delta:=\frac{\ln
M}{\ln \frac{M}{\epsilon}}$$ Since the map $z\mapsto z^N$ is
ergodic, using Birkhoff's theorem we have that, for almost every
$z\in\mathbb{T}$,
$$\lim_{m\rightarrow\infty}\frac{1}{m}\sum_{k=0}^{m-1}\chi_{A_\epsilon}(z^{N^k})=\mu(A_\epsilon).$$
Then, using (\ref{eq4_5_2}), there exists $m_x$ and some
$\delta'>\delta$ such that, for $m\geq m_x$,
$$\frac{1}{m}\sum_{k=0}^{m-1}\chi_{A_\epsilon}(z^{N^k})>\delta'.$$
This can be rewritten as
$$k_{\epsilon,z}:=\sharp\{k\leq
m-1\,|\,z^{N^k}\in A_\epsilon\}>m\delta'.$$ Therefore we have:
\begin{align*}
\prod_{k=0}^{m-1}|m_0(z^{N^k})|^2&\leq\epsilon^{k_{\epsilon,z}}M^{m-k_{\epsilon,z}}\leq\epsilon^{m\delta'}M^{m(1-\delta')}\\
&=e^{m(\delta'\ln\epsilon+(1-\delta')\ln M)}
\end{align*}
Since $\delta'>\delta=\ln M/(\ln(M/\epsilon))$, it follows that
$\gamma:=\delta'\ln\epsilon+(1-\delta')\ln M<0$. Therefore
$$\sum_{m=m_x}^\infty\prod_{k=0}^{m-1}|m_0(z^{N^k})|^2\leq\sum_{m=m_x}^\infty e^{m\gamma}<\infty.$$
Take $z=e^{-2\pi\mbox{i}(x/N^l+\theta_{-l})}$ with $x\in\br$ and
$l(=k_x)\in\mathbb{Z}$, and note that $z^{N^k}=e^{-2\pi
i(x/N^{l-k}+\theta_{-l+k})}$; we get that
$$\sum_{m=m_x}^\infty\left|m_0(\frac{x}{N^{k_x}}+\theta_{-k_x})m_0(\frac{x}{N^{k_x-1}}+\theta_{-k_x+1})...m_0(\frac{x}{N^{k_x-m}}+\theta_{-k_x+m})\right|^2<\infty$$

Neglecting the first few terms, and reversing the order in the products, the equation (\ref{eq4_4_0a}) is
obtained and the corollary is proved.
\par
 We remark that, if $\epsilon=0$, the zero set has  positive measure and the ergodicity implies that almost every trajectory
must go through the zero set, therefore the products appearing in
(\ref{eq4_4_0a}) are zero so we are summing zero terms.

\subsection{\label{frac}Representations on fractals}
\par
We saw in \cite{DutJo} that, if one takes
$$N=3\quad\mbox{and}\quad m_0(z)=\frac{1+z^2}{\sqrt{2}},$$
then the representation associated to this filter (as in theorem
\ref{th0_1}) can be constructed on a Hausdorff measure.
\par
More precisely, consider $\mathcal{R}$, the set of all real
numbers that have a base 3 expansion containing only finitely many
$1$'s. On this set take the Hausdorff measure $\mathcal{H}^s$ with
$s=\log_32$-the Hausdorff dimension of the triadic Cantor set.
Define the unitary operators
$$Uf(x)=\frac{1}{\sqrt{2}}f\left(\frac{x}{3}\right),\quad(x\in\mathcal{R},f\in
L^2(\mathcal{R},\mathcal{H}^s)),$$
$$Tf(x)=f(x-1),\quad(x\in\mathcal{R},f\in
L^2(\mathcal{R},\mathcal{H}^s)),$$ and
$\varphi=\chi_{\mathbf{C}}$, where $\mathbf{C}$ is the triadic
Cantor set.
\par
We will use the techniques developed in the previous sections to
analyze in more detail this representation, give it another form
and then to define a possible Fourier transform related to the
Hausdorff measure $\mathcal{H}^s$.
\par
For our $m_0$, the corresponding function $W=|m_0|^2/3$ is
$$W(z)=\frac13+\frac16z^2+\frac16z^{-2},$$
therefore the coefficients are $a_0=\frac13$, $a_2=a_{-2}=\frac16$
and all others are 0.
\par
We analyze the measure $m$ on $\SN$ constructed from $W$ as in
theorem \ref{th3_3}.

\begin{proposition}\label{prop4_2_1}
The measure $m$ has the following Fourier coefficients:
\begin{equation}\label{eq4_2_1_0}
\hat m(\lambda)=\left\{\begin{array}{ccc}
2^{-(|d_0|+|d_1|+...+|d_p|)/2},&\mbox{if
}\lambda=\sum_{k=0}^p\frac{d_k}{3^k},\mbox{
with }d_k\in\{-2,0,2\},\\
0,&\mbox{otherwise}.
\end{array}
\right.
\end{equation}
\end{proposition}
\begin{proof}
Equation (\ref{eq3_2_4}) becomes in our case
\begin{equation}\label{eq4_2_1_1}
\hat m(\lambda)=\frac12\hat m(3\lambda-2)+\hat
m(3\lambda)+\frac12\hat m(3\lambda+2),\quad(\lambda\in\bz[1/3]).
\end{equation}
We will use another lemma, which is also interesting in its own,
because it tells something about the geometry of the Cantor set
\begin{lemma}\label{lem4_2_2}
$$\hat
m(\lambda)=\mathcal{H}^s((\mathbf{C}+\lambda)\cap\mathbf{C}),\quad(\lambda\in\bz[1/3]).$$
\end{lemma}
\begin{proof}{\it of lemma \ref{lem4_2_2}.} We saw in theorem
\ref{th3_3} that the representation associated to $m_0$ can be
realized on $L^2(m)$ with scaling function $\varphi'=1$. In this
representation, the translation by $\lambda$ is given by
multiplication by the character $\chi_\lambda$, therefore we have
$$\hat m(\lambda)=\ip{T'_\lambda\varphi'}{\varphi'}.$$
However, this representation is isomorphic to the one described
above, the isomorphism mapping the scaling function
$\varphi=\chi_{\mathbf{C}}$ to $\varphi'=1$. In this
representation, translation by $\lambda$ is simply the translation
by $\lambda$ on $\mathcal{R}$. Consequently,
$$\hat
m(\lambda)=\ip{T_\lambda\chi_{\mathbf{C}}}{\chi_{\mathbf{C}}}=\mathcal{H}^s((\mathbf{C}+\lambda)\cap\mathbf{C}).$$
\end{proof}
From lemma \ref{lem4_2_2} we deduce that, if $|\lambda|\geq 1$,
then $\hat m(\lambda)=0$. So if, $|\lambda_1-\lambda_2|\geq 2$
then at most one of them in in the interval $(-1,1)$, so at most
one of them is non-zero. This implies that in the right side of
the equation (\ref{eq4_2_1_1}), at most one of the terms is
non-zero.
\par
Now take $\lambda=l/3^n$. We will proceed by induction on $n$. If
$n=0$, the result follows from proposition \ref{prop3_2}. For
$n>0$, assume that $\hat m(\lambda)$ is not zero. Using
(\ref{eq4_2_1_1}) and the previous statement, we obtain that
exactly one of the terms $\hat m(3\lambda-2),\hat m(3\lambda),\hat
m(3\lambda+2)$ is non-zero. Denote by $\alpha=3\lambda-d_0$ the
one which is not $0$, with $d_0\in\{-2,0,2\}$. Then, analyzing the
three cases, we see that
$$\hat m(\lambda)=\frac{1}{2^{|d_0|/2}}\hat m(\alpha).$$
But $\alpha$ is of the form $k'/3^{n-1}$ so we can use the
induction hypothesis to conclude that $\alpha$ has the form
$\alpha=\sum_{l=1}^pd_l/3^{l-1}$ so
$$\lambda=\sum_{l=0}^p\frac{d_l}{3^l},$$
and (\ref{eq4_2_1_0}) is proved.
\end{proof}

The fact that we have two embodiments of the wavelet
representation associated to the filter
$m_0(z)=\frac{1+z^2}{\sqrt2}$, one on $\mathcal{R}$ with the
Hausdorff measure $\mathcal{H}^s$, and the other on $\SN$ with the
measure $m$, implies that there is an isomorphism between them
which can be interpreted as a Fourier transform on $\mathcal{R}$
since it transforms translations into multiplications:
\begin{corollary}\label{cor4_2_3}
Consider on $\SN$ the measure $m$ that has Fourier coefficients
given in proposition \ref{prop4_2_1}. There is a unique
isomorphism $\mathcal{F}_3$ from $L^2(\mathcal{R},\mathcal{H}^s)$
to $L^2(\SN,m)$ such that for $f\in L^2(\SN,m)$,
$$\mathcal{F}_3T_\lambda\mathcal{F}_3^{-1}f=\chi_{\lambda}f,\quad(\lambda\in\bz[1/3]),$$
$$\mathcal{F}_3U\mathcal{F}_3^{-1}f=m_0f\circ S,$$
$$\mathcal{F}_3\chi_{\mathbf{C}}=1.$$
\end{corollary}

\subsection{\label{sol}$m_0=1$: representations on the solenoid}
We take now $m_0=1$ which obviously satisfies $R_{m_0}1=1$ and we
describe the wavelet representation associated to it. The measure
$m$ associated to $m_0$ as in proposition \ref{prop3_2} verifies
the scaling equations:
$$\hat m(k)=\delta_k,\quad (k\in\bz),\quad \hat m(\lambda)=\hat
m(N\lambda).$$ Therefore
$$\hat m(\lambda)=\delta_\lambda,\quad(\lambda\in\bz[1/N]).$$
But this means that $m$ is the Haar measure $\mu_{\SN}$ on $\SN$.
Hence, from theorem \ref{th3_3} and example \ref{ex2_6}, we have
the following result:
\begin{proposition}\label{prop5_3_1}
Let $H=L^2(\SN,\mu_{\SN})$, $Tf(z_n)_n=z_0f(z_n)_n$, $Uf=f\circ
S$, ($f\in H,(z_n)_n\in\SN$). Let $\varphi$ be the constant
function $1$ on $\SN$. Then $(H,U,T,\varphi)$ is the wavelet
representation associated to $m_0=1$. The representation is
irreducible.
\end{proposition}
With corollary \ref{cor3_6} we have
\begin{corollary}\label{cor5_3_2}
The only functions $h\in\linft$ that satisfy
$$\frac{1}{N}\sum_{k=0}^{N-1}h\left(\frac{\theta+2k\pi}{N}\right)=h(\theta),\quad(\theta\in[-\pi,\pi)),$$
are the constants.
\end{corollary}

\end{document}